\documentclass[11pt]{article}
\usepackage{nicefrac}

\usepackage[utf8]{inputenc}
\usepackage[T1]{fontenc}
\usepackage{lmodern}
\usepackage[english]{babel}
\usepackage{caption}

\usepackage{color}
\usepackage{amsmath}
\catcode`\@=11 \@addtoreset{equation}{section}

\catcode`\@=12

\usepackage{amssymb}
\usepackage{amsfonts}
\usepackage{xcolor}
\usepackage{graphics}
\usepackage{epsfig,psfrag,graphicx}
\usepackage{amssymb}
\usepackage{eepic,epic}
\usepackage{dsfont}

\usepackage{enumerate}

\usepackage{epsfig} 
\usepackage{amsmath} 
\usepackage{amssymb} 
\usepackage{amsbsy} 
\usepackage{comment}

\usepackage{nicefrac}
\usepackage{multirow}
\usepackage{siunitx}
\usepackage[ruled,vlined,linesnumbered]{algorithm2e}

\usepackage{float}

\usepackage{vmargin}
\setmarginsrb  { 0.7in}  
               { 0.6in}  
               { 0.7in}  
               { 0.8in}  
               {  20pt}  
               {0.25in}  
               {   9pt}  
               { 0.3in}  

\DeclareMathOperator*{\argmin}{argmin}

\newtheorem{definition}{Definition}[section]

\newtheorem{lemma}[definition]{Lemma}

\newtheorem{proposition}[definition]{Proposition}
\newtheorem{remark}[definition]{Remark}

\newcommand{\nc}{\newcommand}

\nc{\qed}{\mbox{}\nolinebreak\hfill \rule{2mm}{2mm}}
\nc{\proof}[1]{{\bf Proof of Proposition {#1}: }}

\def\bbbone{{\mathchoice {\rm 1\mskip-4mu l}
{\rm 1\mskip-4mu l} {\rm 1\mskip-4.5mu l} {\rm 1\mskip-5mu l}}}

\def\vec#1{\boldsymbol{#1}}

\newcommand{\dt}		{\, {\rm d}t}

\def\bbbone{{\mathchoice {\rm 1\mskip-4mu l}
{\rm 1\mskip-4mu l} {\rm 1\mskip-4.5mu l} {\rm 1\mskip-5mu l}}}

\renewcommand{\bbbone}{\mathds{1}}

\newcommand{\mc}{\mathcal}

\newcommand{\veps}{\varepsilon}

\newcommand{\oline}{\overline}

\newcommand{\R}{\mathbb{R}}

\newcommand{\N}{\mathbb{N}}

\allowdisplaybreaks

\title{\LARGE \bf{A hybrid level-based learning swarm algorithm with mutation operator for solving large-scale cardinality-constrained portfolio optimization problems
}}
\author{ \textsl{Massimiliano Kaucic}, $\;$\textsl{Filippo Piccotto}, $\;$\textsl{Gabriele Sbaiz}, $\;$\textsl{Giorgio Valentinuz} \vspace{.2cm} \\
\footnotesize{ \textsc{Universit\`a degli Studi di Trieste},} \\
\footnotesize{\textit{Dipartimento di Scienze Economiche, Aziendali, Matematiche e Statistiche “Bruno de Finetti”},} \\
\footnotesize{Via Valerio 4/1, 34127 Trieste, Italy} \vspace{0.1cm} \\
\footnotesize{\ttfamily{massimiliano.kaucic@deams.units.it}$\,,\quad$ 
\ttfamily{filippo.piccotto@phd.units.it}\vspace{.1cm}}\\
\footnotesize{\ttfamily{gabriele.sbaiz@phd.units.it}$\,,\quad$
\ttfamily{giorgio.valentinuz@deams.units.it} \vspace{.1cm}
}
}

\date{\small \today}

\begin{document}
\maketitle

\abstract{
In this work, we propose a hybrid variant of the level-based learning swarm optimizer (LLSO) for solving large-scale portfolio optimization problems.
Our goal is to maximize a modified formulation of the Sharpe ratio subject to cardinality, box and budget constraints. The algorithm involves a projection operator to deal with these three constraints simultaneously and we implicitly control transaction costs thanks to a rebalancing constraint.
We also introduce a suitable exact penalty function to manage the turnover constraint. In addition, we develop an ad hoc mutation operator to modify candidate exemplars in the highest level of the swarm. The experimental results, using three large-scale data sets, show that the inclusion of this procedure improves the accuracy of the solutions. Then, a comparison with other variants of the LLSO algorithm and two state-of-the-art swarm optimizers points out the outstanding performance of the proposed solver in terms of exploration capabilities and solution quality. Finally, we assess the profitability of the portfolio allocation strategy in the last five years using an investible pool of 1119 constituents from the MSCI World Index.
}


\section{Introduction}
In modern portfolio theory, the classical mean-variance portfolio selection problem developed by Markowitz \cite{MW} plays a crucial role. Following this approach, investors should consider together return and risk, distributing the capital among alternative securities based on their return-risk trade-off.
Since the  pioneering work by Markowitz, the mean-variance optimization model (MVO) has been recognized as a practical tool to tackle portfolio optimization problems, and a large number of developments of the basic model has been investigated (see, for instance, \cite{GU} and \cite{KOL}). Several authors have studied the multi-objective formulation of the MVO problem, in which the expected portfolio return is maximized and, at the same time, its variance is minimized (\cite{CE}, \cite{CO} and \cite{KAU}). These contributions aim to provide algorithms able to generate accurate dotted representations of all the sets of non-dominated portfolios in few iterations. On the contrary, in other studies, a single-objective formulation has been considered, maximising the return per unit of risk using a performance measure called the Sharpe ratio (\cite{Sh_1966}, \cite{F-T} and \cite{Z-K}). In this case, the attention is on the selection of the best-performing alternative in the set of feasible portfolios. However, from the point of view of a rational investor, the use of this indicator in periods of market downturns is questionable because it leads to prefer riskier portfolios (\cite{IS} and \cite{K-B-C}).

In this paper, we focus on this occurrence and we propose a maximization problem where the objective function is a modified version of the Sharpe ratio, which is coherent with agent preferences also when risk premia are negative. We further consider four types of real-world constraints. First, a  cardinality constraint manages the size of the portfolio; next, a budget constraint ensures that all the available capital is invested; lastly, bound constraints prescribe lower and upper bounds on the fraction of capital invested in each asset. To conclude, a turnover constraint implicitly controls the effect of the transaction costs on the portfolio rebalancing phases. We analyse the asset allocation problem from the perspective of an institutional investor who operates in large equity markets composed of hundreds or thousands of constituents and selects a restricted pool of stocks to build up a portfolio with a suitable performance with respect to the benchmark.

Since the resulting mixed-integer optimization problem has been proved to be NP-hard, finding possible optimal solutions becomes computationally challenging \cite{E-T-S}. For this reason, on the one hand, exact methods were proposed to supply optimal solutions, but they demand a significant amount of computation time when the problem size increases (\cite{L-S-W}, \cite{S-L-K} and \cite{B-S}). On the other hand, heuristic approaches can identify approximate and sometimes optimal solutions within reasonable computation time even when the problem size is huge (\cite{C-M-B}, \cite{C-S}, \cite{W-L-B} and \cite{L-M}). In this context, swarm optimization algorithms, inspired by the self-organizing interaction among agents, have become popular in portfolio optimization theory in recent years \cite{E-K}. Specifically, the particle swarm optimization (PSO) algorithm has been widely employed to solve real-world financial problems since its first proposal (\cite{CU}, \cite{Z-C-W}, \cite{TH}, \cite{Z-W-W-C} and \cite{C-T-F-P}) due to its effectiveness in reaching optimal solutions. However, the aforementioned algorithm does not work efficiently when the problem size is large, leading to population stagnation and premature convergence \cite{P-V}. To improve PSO performance for large-scale optimization problems, several authors have designed many variants, such as competitive swarm optimizer \cite{C-J_CSO}, social learning particle swarm optimizer \cite{C-J_SLPSO} and level-based learning swarm optimizer (LLSO) \cite{Y-C-D-L-G-Z}. In particular, the latter one has shown better exploitation ability in different environments. For this reason, we propose its use to solve our portfolio optimization problem. LLSO is inspired by the teaching concept that teachers should treat students differently according to their cognitive and learning abilities. Based on that, the general idea of the LLSO is to sort the swarm individuals in ascending order with respect to their fitness and then separate them into distinct levels. The best individuals are stored into the higher level and are not updated, preserving the most valuable information conveyed in the swarm. Unlike PSO, which uses the historically best positions to update the particles, LLSO employs predominant particles in the current swarm to guide the learning of the worst particles and to enhance the swarm diversity. Thus, particles in lower levels have more individuals in the upper levels to learn from and are focused on exploring the search space; those in higher levels mainly concentrate on the exploitation task. Even though LLSO shows promising capabilities in dealing with large-scale optimization problems, it is overly sensitive to its parameters. To mitigate this influence, an adaptive variant, henceforth ALLSO, has been introduced in \cite{S-Y-G-M-L-Z}, which takes advantage of a swarm aggregation indicator to estimate the evolution state of the swarm. Two adaptive adjustment strategies are then applied to identify the best configuration setting for each generation.

Due to the fact that the swarm optimization algorithms are usually blind to the constraints, they have to be equipped with constraint-handling techniques \cite{COE} to be effective in real-world applications. A class of constraint-handling methods widely used in literature is represented by the penalty function methods, where a penalty term reduces the fitness value of the infeasible candidates. However, despite its simplicity, this method usually requires the definition of problem-dependent parameters that significantly impact algorithm performance. To overcome this issue, adaptive penalty techniques have been developed, in which the parameters are automatically set by using  information gathered from the violated constraints at the current generation. For a more exhaustive overview of adaptive penalty techniques, we refer the reader to \cite{B-L}.

Therefore, to tackle the presented large-scale cardinality-constrained portfolio  optimization problem, we combine the ALLSO with a novel hybrid constraint-handling technique, in which we integrate a projection operator into the self-adaptive penalty scheme developed by \cite{C-F-R}. Moreover, to further improve the exploitation power of the algorithm and the quality of the solutions, we introduce a novel mutation procedure, applied to the best individuals in the first level, which generalises the one inspected in \cite{K-B-C}.

Even though similar asset allocation models with real-world constraints have been already studied in the literature, our work represents the first application of the LLSO algorithm for solving a large-scale instance of the problem with the modified Sharpe ratio as the objective function. In addition, to the best of our knowledge, this is the first time such a hybrid constraint-handling technique and the mutation operator are presented and involved in the level-based learning paradigm.

Let us now give a more precise overview of the contents of the paper. The next section describes the investment framework, focusing on the portfolio optimization problem. In Section \ref{sec:opt}, we introduce the developed solver. More precisely, we first explain the adaptive LLSO and then we detail the proposed methods, namely the novel mutation operator and the hybrid constraint-handling technique. In the last part of the section, we summarise the entire procedure. In Section \ref{sec:experiments}, we show the experimental results; in the last section, we depict the conclusions and future perspectives.

\section{Portfolio design}\label{sec:port}


\subsection{Investment framework} \label{subsec:invf}

Let us consider the standard portfolio selection problem introduced in \cite{MW}. We have a frictionless market in which no short selling is allowed, and all investors act as price takers. Assuming that $n$ assets represent the investable universe, a portfolio is identified with the vector of assets weights $\vec x = (x_1,\ldots,x_n)\in \mathbb{R}^n$, where $x_i \in \R$ denotes the proportion of capital invested in asset $i$, with $i = 1,\dots, n$.
Let $R_i$ be the random variable which stands for the rate of return of asset $i$, with expected value $\mu_i$. Hence, the random variable $R_p(\vec x) = \sum_{i=1}^{n}R_i x_i $ indicates the rate of return of portfolio $\vec x$. The expected rate of return of portfolio $\vec x$ is then defined as
\begin{equation}\label{eq:port_return}
  \mu_p(\vec x) = \sum_{i=1}^{n} x_i \mu_i
\end{equation}
and its standard deviation, also called volatility, is given by
\begin{equation}\label{eq:port_risk}
  \sigma_p(\vec x) = \sqrt{\sum_{i=1}^{n} \sum_{j=1}^{n} c_{ij} x_i x_j}
\end{equation}
where $(C)_{ij} = c_{ij}$ is the covariance between stocks $i$ and $j$, with $i, j = 1,\dots, n$.\\Since investors perceive large deviations from the portfolio mean value as damaging, \eqref{eq:port_risk} represents the so-called portfolio risk.

In such a setting, the portfolio choice is made only with respect to the expected portfolio rate of return and the portfolio risk, as stated in the following definition.

\begin{definition}
Given two portfolios $\vec x, \vec y$, we say that $\vec x$ is preferred to $\vec y$ if and only if $\mu_p(\vec x) \geq \mu_p(\vec y)$ and $\sigma_p(\vec x) \leq \sigma_p(\vec y)$, with at least one strict inequality.
\end{definition}

In other words, an investor prefers one portfolio to another if it has a higher expected rate of return and lower risk.


\subsection{Objective function} \label{subsec:objfun}
In our portfolio selection problem, we take the Sharpe ratio (\cite{Sh_1966} and \cite{Sh_1994}) as an essential point of reference for measuring and comparing investment performance. It is defined as the ratio between the expected rate of return of the portfolio in excess of the risk-free rate and the standard deviation of the rates of return of the portfolio itself, that is
\begin{equation}\label{Sharpe}
SR(\vec x) = \frac{\mu_p(\vec x)-r_f}{\sigma_p(\vec x)}
\end{equation}
where $r_f$ is the risk-free rate.
This measure evaluates the compensation earned by the investor per unit of both systematic and idiosyncratic risks \cite{CAP}. Thus, higher values of $SR$ indicate more promising portfolios.

From a theoretical point of view, this choice is justified by the fact that several widely used performance measures are increasing functions of the Sharpe ratio (\cite{S-E_2011} and \cite{S-E_2012}). Moreover, when the numerator in \eqref{Sharpe} is positive, this indicator is coherent with the risk-return profile of a rational investor. From a practical point of view, it can be easily calculated and its interpretation is simpler than most of recently proposed complex performance measures \cite{A-S}.

However, as pointed out in \cite{A-M} and \cite{IS}, the reliability of this performance measure decreases when the excess rate of return is negative. In that case, one would prefer higher-risk portfolios using the Sharpe ratio. To overcome this issue, we adopt the following modification of \eqref{Sharpe}, the so-called modified Sharpe ratio:
\begin{equation}\label{Msr}
MSR(\vec x) = \frac{\mu_p(\vec x)-r_f}{\sigma_p(\vec x)^{\operatorname{sign}(\mu_p(\vec x)-r_f)}}
\end{equation}
where $\text{sign} (z)$ is the sign function of $z \in \mathbb{R}$.
Observe that if the portfolio excess return is non-negative, the modified Sharpe ratio is equal to the Sharpe ratio. Otherwise, it multiplies the portfolio excess return by the standard deviation. In this manner, even in adverse conditions, portfolios with lower risk and higher excess return will be preferred.

\subsection{Constraints} \label{subsec:constrains}
In our portfolio model, we consider the following constraints.
\begin{itemize}
  \item \emph{Budget}. All the available capital needs to be invested. In terms of portfolio weights, this translates to
      \begin{equation}\label{budg}
        \sum_{i=1}^{n} x_i = 1.
      \end{equation}
  \item \emph{Cardinality}. We assume that the portfolio includes up to $k$ assets, where $k\leq n$. To model the inclusion or the exclusion of the $i$-th asset in the portfolio, a binary variable $\delta_i$ is introduced as
      \begin{equation}
        \delta_i =
        \left\{
        \begin{array}{l}
        0\mbox{, if asset $i$ is excluded} \\
        1\mbox{, if asset $i$ is included}
        \end{array}
        \right.
      \end{equation}
      for $i = 1,\ldots,n$. The resulting vector of selected assets is $\vec \delta = (\delta_1,\ldots,\delta_n) \in \{0, 1\}^n$, and the cardinality constraint can be written as
      \begin{equation}\label{card}
      \sum_{i=1}^{n} \delta_i \leq k.
      \end{equation}
  \item \emph{Box}. A balanced portfolio should avoid extreme positions and foster diversification. Hence, we impose a maximum and a minimum limit for portfolio weights, that is
      \begin{equation}\label{box}
        \delta_i l_i \leq x_i \leq \delta_i u_i,\quad i = 1,\ldots,n
      \end{equation}
      where $l_i$ and $u_i$ are the lower and the upper bounds for the weight of the $i$-th asset, respectively, with $0< l_i < u_i \leq 1$ to exclude short sales.
  \item \emph{Turnover}. To control the effect of the transaction costs in the portfolio rebalancing phases, we consider a portfolio turnover constraint. Let $\vec x_0$ be a vector containing the current portfolio positions \cite{S-W-M}. Then, the portfolio turnover constraint is
      \begin{equation}\label{turn}
        \sum_{i=1}^{n} \lvert x_{i} - x_{0 ,\, i} \rvert \leq TR
      \end{equation}
      where $TR$ denotes the maximum turnover rate, which lies between $0$ and $1$. Note that if $TR=0$ rebalancing is not allowed, and more trades are allowed when $TR$ increases.
\end{itemize}
The pairs $(\vec \delta ,\, \vec x) \in \{0, 1\}^n \times \mathbb{R}^n$ that satisfy \eqref{budg}, \eqref{card}, \eqref{box} and \eqref{turn} form the feasible set $\mathcal{F}$. Then, our portfolio optimization problem can be written as
\begin{equation}\label{opt}
\begin{aligned}
\max_{\vec \delta ,\, \vec x} \quad & MSR(\vec x)\\
\textrm{s.t.} \quad & \left(\vec \delta, \vec x\right) \in \mathcal{F}.
\end{aligned}
\end{equation}

\begin{remark}
To simplify the following treatment, we reformulate our maximization problem into the equivalent minimization problem
\begin{equation}\label{opt1}
\begin{aligned}
\min_{\vec \delta, \vec x} \quad & f(\vec x)\\
\textrm{s.t.} \quad & \left(\vec \delta, \vec x\right) \in \mathcal{F}
\end{aligned}
\end{equation}
where $f(\vec x) = -MSR(\vec x)$.
\end{remark}

\section{Optimization algorithm} \label{sec:opt}

\subsection{Adaptive level-based learning swarm optimizer} \label{subsec:allso}

The algorithm evolves a swarm of $NP$ candidate solutions using the so-called level-based population structure \cite{Y-C-D-L-G-Z}, according to which the evolution process is defined as follows.
\begin{enumerate}
\item At each iteration $g$, the individuals in the swarm are first sorted ascending based on their fitness and grouped into $NL_g$ levels, each one containing $LP_g = \lfloor\nicefrac{NP}{NL_g}\rfloor$ particles. In the last level, there are $\lfloor\nicefrac{NP}{NL_g}\rfloor + NP \% NL_g$ particles.\footnote{We denote by $\lfloor x \rfloor$ the floor of $x$ and by $x\%y$ the rest of the division of $x$ by $y$.} Better individuals belong to higher levels, and a higher level corresponds to a smaller level index. Thus, $L_{1}$ represents the best level, while $L_{NL_g}$ is the worst one.
\item To preserve the most valuable information conveyed in the current swarm, individuals belonging to $L_1$ are not updated and enter directly in the next generation. The $p$-th particle in level $L_{l}$, denoted by $\vec x^{l ,\, p}(g)$, where $l = 3,\ldots,NL_g$ and $p = 1,\ldots,LP_g$, is allowed to learn from two particles $\vec x^{l_1 ,\, p_1}(g)$, $\vec x^{l_2 ,\, p_2}(g)$ randomly extracted from two different higher levels $L_{l_1}$ and $L_{l_2}$ with $l_1 < l_2$, and $p_1$ and $p_2$ are randomly chosen from $\{1,\ldots,LP_g\}$. For $l=2$, we sample two particles from $L_{1}$ in such a way that $\vec x^{l_1 ,\, p_1}(g)$ is better than $\vec x^{l_1 ,\, p_2}(g)$ in terms of fitness function. Thus, the update rule for particle $\vec x^{l ,\, p}(g)$ is given component-wise by
    \begin{align}
    &v_{i}^{l ,\, p}(g+1) = r_1 v_{i}^{l ,\, p}(g) + r_2
    \left(x_{i}^{l_1 ,\, p_1}(g)-x_{i}^{l ,\, p}(g)\right) + \phi_g r_3  \left(x_{i}^{l_2 ,\, p_2}(g)-x_{i}^{l ,\, p}(g)\right) \label{update:v}\\
    &x_{i}^{l ,\, p}(g+1) = x_{i}^{l ,\, p}(g) + v_{i}^{l ,\, p}(g) \label{update:x}
    \end{align}
    for $i = 1, \ldots, n$, where $v_{i}^{l ,\, p}(g)$ denotes the $i$-th component of the velocity of particle $p$ in level $L_l$ at generation $g$, and $r_1$, $r_2$, $r_3$ are real numbers randomly generated within $[0,1]$. The parameter $\phi_g \in [0,1]$ controls the influence of the less performing exemplar $\vec x^{l_2 ,\, p_2}(g)$ on $\vec v^{l ,\, p}(g)$.
\end{enumerate}
Based on \cite{S-Y-G-M-L-Z}, both the parameters involved in the learning process at generation $g$, namely $NL_g$ and $\phi_g$, are adaptively adjusted based on the evolution state of the swarm by an aggregation indicator, which is defined as
\begin{equation} \label{EQN:aggregation_indicator}
s(g) = \frac{\oline f_g-f\left(\vec x_{gbest}(g)\right)}{f\left(\vec x_{gbest}(g)\right)+ \xi}
\end{equation}
where $\oline f_g$ is the average fitness of the population at generation $g$, $f\left(\vec x_{gbest}(g)\right)$ denotes the historically global best fitness up to iteration $g$, and $\xi$ is a small positive value to avoid zero denominators.

\begin{remark}
When $s(g)$ is high, particles are far from the current global best solution. Thus, the swarm is in an exploration phase. On the contrary, when $s(g)$ is low, particles are close to the global best solution $\vec x_{gbest}(g)$ and the swarm is in an exploitation phase.
\end{remark}

To guarantee a control on the number of levels, $NL_g$ takes values in the set $\{NL_{min}, \ldots, NL_{max}\}$, where $NL_{min} ,\, NL_{max} \in \mathbb{N}$ are predefined lower and upper bounds. Moreover, to balance the level selection diversity and the exemplar diversity, $NL_g$ can be modified only when the relative improvement of the global fitness between generation $g$ and generation $g-1$, given by
\begin{equation} \label{eqn:relative_improvement}
t(g) = \frac{f\left(\vec x_{gbest}(g-1)\right)-f\left(\vec x_{gbest}(g)\right)}{f\left(\vec x_{gbest}(g)\right)+\xi},
\end{equation}
slows down or stops, which corresponds to the cases $t(g) < t(g-1)$ or $t(g) = 0$ respectively. The update of $NL_g$ then follows the rule
\begin{equation} \label{update:NL1}
{NL}_g = \begin{cases}
2\cdot {NL}_{g-1} \quad \text{ if }\quad s(g) < \bar{\delta} \\
\frac{1}{2}\cdot {NL}_{g-1} \quad \text{ if }\quad s(g) \geq \bar{\delta}
\end{cases}
\end{equation}
where $\bar{\delta}$ is a threshold in terms of the aggregation indicator to control the adjustment of $NL_g$.

When $NL_g$ is out of the range, it is adjusted as follows
\begin{equation} \label{update:NL2}
{NL}_g = \begin{cases} {NL}_{rand} \quad \text{ if }\quad r<px\\
NL_{max} \quad \text{ if }\quad r \geq px \text{ and } {NL}_g>{NL}_{max}\\
NL_{min} \quad \text{ if }\quad r \geq px \text{ and } {NL}_g<{NL}_{min}
\end{cases}
\end{equation}
where $NL_{rand}$ is uniformly sampled from $\{NL_{min}, \ldots, NL_{max}\}$, $r$ is a real number randomly generated within $[0,1]$, and $px$ is a fixed probability employed to reset $NL_g$.

The update for $\phi_g$ is designed in the following way
\begin{equation} \label{update:phi_g}
\phi_{g} = 0.35+0.1\cdot \frac{1}{1+10\cdot s(g)}
\end{equation}
where $s(g)$ is the value of aggregation indicator given in \eqref{EQN:aggregation_indicator}.

A preliminary numerical analysis reveals that a clamping procedure, limiting the magnitude of the velocity $\vec v^{l ,\, p}(g)$, provide a better exploration of the search space (in this regard, see also \cite{O-E-C}). This function can be written component-wise as
\begin{equation}\label{cla}
v^{l ,\, p}_{i}(g) = \min \{\max\{v^{l ,\, p}_{i}(g),v_i^{min}\},v_i^{max}\},
\end{equation}
where $v_i^{min}$ and $v_i^{max}$ are the minimum and the maximum velocity allowed for component $i$, with $i = 1, \ldots, n$. In the experiments, recalling equation \eqref{box}, we set $v_i^{max} = u_i$ and $v_i^{min} = -v_i^{max}$.

\subsection{Mutation operator}\label{subsec:mut}
Instead of directly moving the individuals of the first level to the next generation, we propose to mutate them using an operator that combines two perturbation strategies properly developed for our portfolio optimization problem.\\More specifically, one technique is inspired by the swap operator proposed in \cite{K-M-P} and works as follows. First, we fix the maximum allowed number of non-null positions that could become zero, namely $k_{max}^{swap}$. Then, for each particle $\vec x^{1 ,\, p}(g)$ in level $L_1$ subject to swapping, we randomly sample from $\left\{1, \ldots, k_{max}^{swap}\right\}$ the number $k^{\textrm swap}$ of non-null positions that will be set to zero. At this point, for $j = 1, \ldots, k^{swap}$, let $a_j$ and $b_j$ be two randomly chosen positions in $\vec x^{1 ,\, p}(g)$, such that $x_{a_j}^{1 ,\, p}(g) = 0$ and $l_{b_j} \leq x_{b_j}^{1 ,\, p}(g) \leq u_{b_j}$. Thus, the modified individual, $\hat{\vec x}^{1 ,\, p}(g)$, is defined component-wise as
\begin{equation}\label{PSO:swap_operator}
    \hat{x}_{i}^{1 ,\, p}(g) = \left\{
      \begin{array}{ll}
        x_{i}^{1 ,\, p}(g), & \hbox{if $i \neq a_j \text{\ and\ } i \neq b_j$} \\
        l_{a_j} + \dfrac{x_{b_j}^{1 ,\, p}(g)-l_{b_j}}{u_{b_j}-l_{b_j}}(u_{a_j}-l_{a_j}), & \hbox{if $i = a_j$} \\
        0, & \hbox{if $i = b_j$}.
      \end{array}
    \right.
\end{equation}
In this paper, based on the preliminary experiments, $k_{max}^{swap} = \lfloor 0.05 \cdot k\rfloor$, where $k$ represents the maximum number of assets included in the portfolio.

\begin{remark}
This generalisation, allowing multiple swaps at the same time, improves the search capabilities of the original swap operator.
\end{remark}

The other perturbation scheme focuses solely on the non-null components. For each $\vec x^{1 ,\, p}(g)$ to be mutated, let $I_+^{1 ,\, p}(g) = \{i : x_i^{1 ,\, p}(g) > 0\}$ then, for all $i \in I_+^{1 ,\, p}(g)$, we define the interval
\begin{equation}
  W_i^{1 ,\, p}(g) = \left[x_i^{1 ,\, p}(g) - \Delta_i(g) ,\, x_i^{1 ,\, p}(g) + \Delta_i(g)\right]
\end{equation}
where $\Delta(g) = \left(1-\frac{g}{g_{max}+1}\right) \left(\mathbf{u} - \mathbf{l}\right)$, with $g_{max}$ be the maximum allowed number of iterations. The mutated component $\hat{x}_i^{1 ,\, p}(g)$ is randomly generated from the interval $W_i^{1 ,\, p}(g) \cap \left[l_i, u_i\right]$. For $i \notin I_+^{1 ,\, p}(g)$, we set $\hat{x}_i^{1 ,\, p}(g) = 0$. By narrowing the range over time, this procedure increases the exploration around the particles in $L_1$.

\begin{remark}
By construction, the solutions modified by both the perturbation operators have at most $k$ non-null positions and satisfy the box constraints.
\end{remark}

For each particle in $L_1$, the probability of applying the generalised swap operator decreases as the iteration counter increases according to the following rule
\begin{equation}\label{EQN:pSwap}
p_{swap}(g) = \frac{1}{1+\exp(-0.005 \cdot g)}.
\end{equation}
In the initial stages, the proposed mutation favours the global search, using the generalised swap operator to identify the most promising subset of non-null decision variables. With the progress of the generations, the role of the refinement operator increases and, in the late stages, the algorithm focuses primarily on the local search.

The pseudo-code of the developed mutation procedure is reported in Algorithm \ref{PSEUDOCODE:local_search}.

\resizebox{0.85\textwidth}{!}{
\begin{algorithm}[H]
\DontPrintSemicolon
\SetAlgoLined
\SetKwInOut{Input}{Input}\SetKwInOut{Output}{Output}
\Input{$\vec x^{1 ,\, p}(g)$, $\vec l$, $\vec u$, $k_{max}^{swap}$, $g$, $\Delta(g)$}
\Output{$\hat{\vec x}^{1 ,\, p}(g)$}
\BlankLine
Set $\hat{\vec x}^{1 ,\, p}(g) = \vec x^{1 ,\, p}(g)$ \;
Set $I^0 = \left\{i \colon x_i^{1 ,\, p}(g) = 0\right\}$ \;
Set $I^+ = \left\{i \colon l_i \leq x_i^{1 ,\, p}(g) \leq u_i \right\}$\;
Calculate $p_{\textrm swap}(g)$ according to \eqref{EQN:pSwap} \;
\uIf{$rand() \leq p_{\textrm swap}(g)$}{
    $k^{swap} \rightarrow \left\{1, \ldots, k_{max}^{swap}\right\}$ \;
    \For{$j = 1$ to $k^{swap}$}{
    $a_j \rightarrow I^0$ \;
    $b_j \rightarrow I^+$ \;
    $\hat{x}^{1 ,\, p}_{a_j}(g) = l_{a_j} + \dfrac{x_{b_j}^{1 ,\, p}(g)-l_{b_j}}{u_{b_j}-l_{b_j}}(u_{a_j}-l_{a_j})$ \;
    $\hat{x}^{1 ,\, p}_{b_j}(g) = 0$ \;
    }
}
\Else{
    \For{$i$ in $I^+$}{
        $lb = \max\left(x^{1 ,\, p}_i(g) - \Delta_i(g), l_i\right)$ \;
        $ub = \min\left(x^{1 ,\, p}_i(g) + \Delta_i(g), u_i\right)$ \;
        $\hat{x}^{1 ,\, p}_i(g) \rightarrow [lb, ub]$ \;
    }
}
\caption{Mutation procedure} \label{PSEUDOCODE:local_search}
\end{algorithm}
}

\subsection{Solution coding and hybrid constraint-handling procedure} \label{subsec:consthndling}

Let us introduce some notation. Let $\mc C_i$ denote a closed convex subset of $\R_+$, with $i = 1, \ldots, n$, and $K = \left\{i_1, \ldots, i_k\right\}$ be any subset of indices of  $I=\{1, \ldots, n\}$ with cardinality $k \in \N$, so that $I \setminus K$ is its complement in $I$. For all $\vec x \in \R^n$, let $\vec x_K$ be defined component-wise as
\begin{equation} \label{coding:real}
\vec x_{K ,\, i} =
\left\{
  \begin{array}{ll}
    x_i, & \hbox{if $i \in K$} \\
    0, & \hbox{if $i \in I \setminus K$}
  \end{array}
\right.
\end{equation}
and let $\pi_K \colon \R^n \to \R^k$ be the projection such that $\pi_K(\vec x) = \left(x_{i_1}, \ldots, x_{i_k}\right)$.

We start by presenting the following proposition (in this regard, see also \cite{Z-L-A}).
\begin{proposition} \label{prop_proj}
Let $\vec y \in \R^n$, with $n\geq 2$. Then, the optimal $K$ for the problem
\begin{equation} \label{probl_components}
\min_{\vec x_K:\, x_i \in C_i} \frac{1}{2}\|\vec x_K-\vec y\|^2
\end{equation}
is the set $K^\ast$ of indices corresponding to the $k$ largest components of $\vec y$.
\end{proposition}
The proof of this result is reported in Appendix \ref{Appendix_A}.\\
In other words, the proposition states that $\vec x_{K^\ast}$ is the vector with at most $k$ non-null components which has minimum Euclidean distance from $\vec y$ among all $\vec x_K$, with $K \subset I$ of cardinality $k$.\\
Thanks to this projection, which implicitly enforces cardinality fulfillment, we can remove the vector of binary variables $\vec \delta$ from the coding scheme of the solutions and we reformulate the portfolio optimization problem \eqref{opt1} only in terms of $\vec x_{K^\ast}$.
To this end, we introduce the set
\begin{equation} \label{feasible_set}
\mc B = \left\{\vec x \in \R^n \colon x_i = 0 \text{ or } x_i \in \left[l_i ,\, u_i\right] \text{ for } i \in K^\ast ,\, x_i = 0 \text{ for } i \in I \setminus K^\ast ,\,  \sum_{i = 1}^n x_i = 1\right\}
\end{equation}
that is the set of the points satisfying all the constraints apart from the turnover condition. Further, let $\psi(\vec x)$ represent the value of the turnover function at $\vec x$, which is given by
\begin{equation}\label{turn_value}
\psi(\vec x) = \sum_{i=1}^{n} \lvert x_i - x_{0,i} \rvert - TR.
\end{equation}
Then, the constrained optimization problem can be rewritten as
\begin{equation}\label{opt_bis}
\begin{aligned}
\min_{\vec x \in \mc B} \quad & f(\vec x)\\
\textrm{s.t.} \quad & \psi(\vec x) \leq 0.
\end{aligned}
\end{equation}
The following proposition, whose proof is given in Appendix \ref{Appendix_A}, establishes the equivalence between problem \eqref{opt_bis} and the mixed-integer optimization problem \eqref{opt1}.
\begin{proposition} \label{prop_equiv}
We assume that $\left(\vec \delta^\ast ,\, \vec x^\ast\right)$ is a global solution to problem \eqref{opt1}, then $\vec x_{K^\ast}^\ast$ is a global solution to problem \eqref{opt_bis}. Conversely, if $\vec x^\ast$ is a global solution to \eqref{opt_bis}, then $\left(\vec \delta^\ast ,\, \vec x^\ast\right)$ is a global solution to \eqref{opt1}, with
\begin{equation*}
\delta_i^\ast =\begin{cases}
1\quad \text{if} \quad i\in K^\ast\\
0 \quad \text{otherwise}.
\end{cases}
\end{equation*}
\end{proposition}

\bigskip

As previously observed, the standard ALLSO algorithm can only deal with unconstrained problems; thus, we propose to incorporate a hybrid constraint-handling technique in order to solve problem \eqref{opt_bis}.

The building block of our procedure is based on the following lemma.
\begin{lemma}
Let $\vec l = \left(l_1, \ldots, l_n\right)$ and $\vec u = \left(u_1, \ldots, u_n\right)$ be such that $l_i \leq u_i$ for $i = 1, \ldots, n$. Let $\vec y \in \R^n$ and define $[\vec l ,\, \vec u] = \left\{\vec x \in \R^n \colon l_i \leq x_i \leq u_i\right\}$. Then, the orthogonal projection of $\vec y$ onto $[\vec l ,\, \vec u]$ is given component-wise by
\begin{equation}\label{proj:box}
P_{[\vec l ,\, \vec u] ,\, i}(\vec y) = \min\{\max\{y_i ,\, l_i\} ,\, u_i\}
\end{equation}
with $i = 1, \ldots, n$.
\end{lemma}
The derivation of the orthogonal projection $P_{[\vec l ,\, \vec u]}$ can be found in \cite{B}. We now provide the main result concerning the projection phase. We refer the reader to Appendix \ref{Appendix_A} for the proof.

\begin{proposition}\label{prop:proj_beck}
Let $\vec y \in \R^n$, with $n \geq 2$, and $K^{\ast\ast} = \left\{i \in K^{\ast} \colon y_i > 0\right\}$, with $K^{\ast}$ being the optimal set in Proposition \ref{prop_proj}. Assume that $\mc B$ in \eqref{feasible_set} is non-empty. Then, the orthogonal projection of $\vec y$ onto $\mc B$ is
\begin{equation}\label{eqn:proj_onto_B1}
P_{\mc B}(\vec y) = \pi_{K^{\ast\ast}}^{-1}\left(P_{[ \pi_{K^{\ast\ast}}(\vec l) ,\,  \pi_{K^{\ast\ast}}(\vec u)]}( \pi_{K^{\ast\ast}}(\vec y - \eta^{\ast} \vec 1)\right)
\end{equation}
where $\pi_{K}^{-1}(\vec z)$ is the pre-image of $\vec z \in \R^{\left|K^{\ast\ast}\right|}$ under $\pi_{K}$ and $\eta^{\ast} \in \R$ is a solution of
\begin{equation}\label{eqn:proj_onto_B2}
\sum_{i=1}^k  P_{[ \pi_{K^{\ast\ast}}(\vec l) ,\,  \pi_{K^{\ast\ast}}(\vec u)]}(\vec y - \eta \vec 1 )=1.
\end{equation}
\end{proposition}
Let $\mc P_g = \left\{\vec x^p(g)\in \R^n \colon p = 1, \ldots, NP\right\}$ be the swarm at generation $g$, with $g = 1, \ldots, g_{max}$. Then, the proposed ALLSO variant maps the individuals in $\mc P_g$, which are updated using \eqref{update:v} and \eqref{update:x}, onto the set $\mc B$ by means of the projector defined in \eqref{eqn:proj_onto_B1}. The resulting mutated swarm is denoted by $\check{\mc P}_g$. Successively, we apply the self-adaptive penalty approach by \cite{C-F-R} to handle the turnover constraint and to guarantee the global optimality of solutions.
More precisely, the objective function value at each projected individual in $\check{\mc P}_g$, namely $f(\check{\vec x}^p)$, is normalized according to the formula
\begin{equation*}
\hat{f}(\check{\vec x}^p) = \frac{f(\check{\vec x}^p)-f^{min}}{f^{max}-f^{min}}
\end{equation*}
where $f^{min} = \displaystyle \min_{\check{\vec x}^p \in \check{\mc P}_g} f(\check{\vec x}^p)$ and $f^{max} = \displaystyle \max_{\check{\vec x}^p \in \check{\mc P}_g} f(\check{\vec x}^p)$. Similarly,
the corresponding normalized constraint violation is given by
\begin{equation*}
  \Psi\left(\check{\vec x}^p\right) =
\left\{
  \begin{array}{ll}
    \dfrac{\max \{\psi(\check{\vec x}^p), 0\}}{\psi^{max}}, & \hbox{if $\psi^{max} > 0$} \\
    0, & \hbox{otherwise}
  \end{array}
\right.
\end{equation*}
where $\psi^{max}$ denotes the maximum of $\psi(\check{\vec x}^p)$ over all the mutated solutions in $\check{\mc P}_g$ which do not satisfy the turnover constraint.

Finally, the penalty function is
\begin{equation}\label{penalty_fun}
F\left(\check{\vec x}^p\right)=\begin{cases}
\hat{f}\left(\check{\vec x}^p\right)\qquad \qquad \quad \; \; \text{if}\quad \psi(\check{\vec x}^p) \leq 0 \\
\hat{f}(\check{\vec z}) + R_f \Psi\left(\check{\vec x}^p\right) \quad \text{if}\quad  \psi(\check{\vec x}^p) > 0 \; \text{ and }\; f(\check{\vec x}^p)\leq f(\check{\vec z})\\
\hat{f}(\check{\vec x}^p) + R_f \Psi\left(\check{\vec x}^p\right) \quad \text{if}\quad \psi(\check{\vec x}^p) > 0 \; \text{ and }\; f(\check{\vec x}^p)> f(\check{\vec z})\, ,
\end{cases}
\end{equation}
where $R_f$ represents the feasibility ratio for $\check{\mc P}_g$, that is the percentage of individuals in $\check{\mc P}_g$ satisfying the turnover constraint. In \eqref{penalty_fun}, the reference point $\check{\vec z}$ is a point belonging to $\check{\mc P}_g$ that satisfies the turnover constraint and has the lowest objective function value found so far. As in \cite{C-F-R}, if the population has no feasible points, $f\left(\check{\vec z}\right)$ is initially and temporarily set to $f^{max}$, so that $f(\check{\vec x}^p)\leq f(\check{\vec z})$ for all $\check{\vec x}^p \in \check{\mc P}_g$ and $\hat{f}(\check{\vec z}) = 1$. The value of $f\left(\check{\vec z}\right)$ is updated only when the first feasible point is encountered.

We conclude this subsection by stating the following theorem, whose proof is omitted since it is similar to the one presented in \cite{C-F-R}.
\begin{proposition}
The problem
$$ \min_{\vec x\in \mc B}F(\vec x)$$
with $F$ as in \eqref{penalty_fun}, is equivalent to the problem \eqref{opt_bis}.
\end{proposition}

\subsection{Initialisation strategy and complete algorithm} \label{subsec:pseudo-code}

For high dimensional problems, the common strategies of seeking a search space coverage by initializing the particles uniformly throughout the space as well as by increasing the size of the swarm are inefficient, because the search space grows exponentially with the dimension \cite{V-E}. Moreover, the presence of highly constrained feasible regions in the search space exacerbates even more the initialization issue \cite{D-D}.

To effectively address the low degree of feasibility in our portfolio rebalancing problem due to the complexity of the turnover constraint, we propose a direct initialization of the candidate solutions in a neighbourhood of $\vec x_0$.\\
Let $d^{min}_i$ and $d^{max}_i$ be the minimum and the maximum allowed weight changes for $\vec x_{0 ,\, i}$ respectively, with $i = 1, \ldots, n$. Let $D^p$ denote the total portfolio weight allowed to be re-allocated in $\vec x_0$ for defining the $p$-th candidate solution $\vec x^p(0)$, with $p = 1, \ldots, NP$. Then, for each $p$,
\begin{enumerate}
  \item we randomly select $D^p$ within $\left[0 ,\, \nicefrac{TR}{2}\right]$;
  \item we select a subset $J^-$ of $k'$ assets from the $k$ assets with positive weight in $\vec x_0$, so that
      $$x^p_j(0) = x_{0 ,\, j} - d_j ,\, \text{ for } j \in J^-$$
      where $d_j$ is randomly sampled in $\left[d^{min}_j ,\, d^{max}_j\right]$ in such a way that $\sum_{j \in J^-} d_j = D^p$, and $x^p_j(0) = 0$ or $l_j \leq x^p_j(0) \leq u_j$;
  \item we select a subset $J^+$ of $k''$ assets from the $n-k$ assets with zero weight in $\vec x_0$, with $k'' \leq k'$, so that
      $$x^p_j(0) = x_{0 ,\, j} + d_j ,\, \text{ for } j \in J^+$$
      where $d_j$ is randomly sampled in $\left[d^{min} ,\, d^{max}\right]$ in such a way that $\sum_{j \in J^+} d_j = D^p$, and $l_j \leq x_j \leq u_j$;
  \item for $j \in I \setminus \left(J^- \cup J^+\right)$, we set $x^p_j(0) = x_{0 ,\, j}$.
\end{enumerate}
The portfolios assembled using this scheme satisfy cardinality, box and turnover constraints. In this way, the initialization strategy encourages the swarm to focus on exploitation rather than exploration, thereby allowing it to identify promising solutions, even in problems with high dimension and small feasible regions.\\Regarding the initial velocities, we set them all equal to the zero vector, that is $\vec v^p(0) = \vec 0$, for $p = 1, \ldots, NP$.

The pseudocode of the proposed LLSO variant with adaptive parameters update, mutation of the particles in the first level and hybrid constraint-handling technique, shortly ALLSO-MUT-H, is reported in Algorithm \ref{PSEUDOCODE:ALLSO-MUT-H}. It can be noticed that, setting $\vec x_0 = \vec 0$ and $TR = 1$, ALLSO-MUT-H can also tackle portfolio optimization problems with no rebalancing. In this case, only the orthogonal projector is needed to move the unfeasible solutions to the feasible region.

\resizebox{0.8\textwidth}{!}{
\begin{algorithm}[H]
\DontPrintSemicolon
\SetAlgoLined
\SetKwInOut{Input}{Input}\SetKwInOut{Output}{Output}
\Input{$\vec x_0$, $\vec l$, $\vec u$, $k$, $TR$, $NL_{min}$, $NL_{max}$, $NP$, $\bar{\delta}$, $px$, $\xi$}
\Output{$\vec x_{gbest}$}
\BlankLine
Set $g = 0$ and $NL_g = 20$ \;
Initialize the swarm $\mc P_g = \left\{\vec x^p(g) \colon p = 1, \ldots, NP\right\}$ and the velocities $\vec v^p(g)$ \;
\For{$i = 1$ to $NP$}{
    Project $\vec x^p(g)$ onto $\mc B$ using \eqref{eqn:proj_onto_B1} \;
    Calculate the turnover violation using \eqref{turn_value} \;
}
Calculate the penalty $F$ for particles in $\mc P_g$ \;
Sort $\mc P_g$ by $F$ value and divide it in $NL_g$ levels \;
Set $\vec x_{gbest}(g) = \vec x^1(g)$ \;
\While{$g < g_{max}$}{
    $g = g + 1$ \;
    Set $\mc P^{L1} = \left\{\vec x^p(g) \colon p \in L_1\right\}$ \;
    \For{$p = 1$ to $LP_{g}$}{
        Use Algorithm \ref{PSEUDOCODE:local_search} to generate the mutated particle $\hat{\vec x}^{1 ,\, p}(g)$ from $\vec x^{1 ,\, p}(g)$ \;
        Project $\hat{\vec x}^{1 ,\, p}(g)$ onto $\mc B$ using \eqref{eqn:proj_onto_B1} \;
        Calculate the turnover violation using \eqref{turn_value} \;
    }
    Set $\mc P^{L1}_{mut} = \left\{\hat{\vec x}^p(g) \colon p \in L_1\right\}$  \;
    Calculate the penalty $F$ for $\mc P^{L1} \cup \mc P^{L1}_{mut}$ and update $\mc P^{L1}$ based on $\mc P^{L1}_{mut}$ using $F$ \;
    Sort $\mc P^{L1}$ by the $F$ value \;
    Calculate the swarm aggregation indicator using \eqref{EQN:aggregation_indicator} and update $\phi_g$ using \eqref{update:phi_g} \;
    \For{$p = LP_{g+1}$ to $NP$}{
        Update $\vec v^p(g)$ using \eqref{update:v} and clamp it using \eqref{cla} \;
        Update $\vec x^p(g)$ using \eqref{update:x} \;
        Project $\vec x^p(g)$ onto $\mc B$ \;
        Calculate the turnover violation using \eqref{turn_value} \;
    }
    Set $\check{\mc P}_g$ be the set of updated particles \;
    Calculate $\phi$ for $\mc P_g \cup \check{\mc P}_g$ and update $\mc P_g$ based on $\check{\mc P}_g$ using $F$ \;
    Sort $\mc P_g$ by $F$ value \;
    Calculate $F$ for the set $\left\{\vec x_{gbest}(g), \vec x^{1}(g)\right\}$ and update $\vec x_{gbest}(g)$ \;
    Calculate the relative improvement $t(g)$ of $\vec x_{gbest}(g)$ using \eqref{eqn:relative_improvement} \;
    \uIf{$t(g)<t(g-1)$ or $t(g)=0$}{
        Update $NL_g$ using \eqref{update:NL1} and \eqref{update:NL2} \;
    }
}
\caption{ALLSO-MUT-H} \label{PSEUDOCODE:ALLSO-MUT-H}
\end{algorithm}
}

\section{Experimental analysis}\label{sec:experiments}
This section is divided into two parts. We first point out the strengths and weaknesses of using the proposed algorithm to tackle large-scale cardinality-constrained portfolio optimization problems. The comparisons are made with two recent variants of the LLSO as well as with other state-of-the-art swarm optimization algorithms implementing the exact $\ell_1$-penalty function approach proposed in \cite{C-T-F-P}. Finally, we assess the profitability of the investment strategy in a real-world case study by varying the size of portfolios.
\subsection{Algorithmic comparisons} \label{subsec:algo_comparisons}
For the algorithmic comparisons, we use three data sets from the OR-Library \cite{C-B}, namely S$\&$P $500$, Russell 2000 and Russell 3000, which represent large capital market indices. 
In Table \ref{Table:1} are summarized the above cited data sets. For the calculation of the expected rates of return, we adopt a historical approach based on all the information available, that consists of 290 weekly prices for each asset. Then, to reduce the bias in the estimation of the covariance matrix $C$, we use the shrinkage estimator proposed in \cite{L-W}.
\begin{table}[ht]
\begin{center}
\captionsetup{width=.9\linewidth}
\caption{Data sets from \cite{C-B} with the corresponding number of weeks and number of market constituents ($n$) used in the estimation of parameters.} \label{Table:1}
{\small
\begin{tabular}{lccccccc}
\hline
Data set name & Weekly prices & Assets ($n$) \\
\hline
S\&P 500 & 290 & 457 \\
Russell 2000 & 290 & 1318 \\
Russell 3000 & 290 & 2151 \\
\hline
\end{tabular}
}
\end{center}
\end{table}

For a fair comparison of the solvers, all the algorithms have the same initial population of $500$ individuals for each test set. In particular, we set $d^{min}_i = 0.0005$ and $d^{max}_i = 0.0050$ in the initialization strategy. For each algorithm, we perform $25$ independent runs with $2000$ generations.
Moreover, all the portfolios from a given test set employ the following parameter setting. The risk-free value in \eqref{Msr} is set to zero. For the cardinality constraint \eqref{card}, we consider $k$ equal to $30\%$ of the size $n$ of the corresponding data set. The box thresholds in \eqref{box} are $l_i=0.001$ and $u_i=0.05$ for each asset. Regarding the turnover constraint \eqref{turn}, $TR$ is set equal to $0.20$ and the vector of current positions $\vec x_0$ is fixed for all the compared algorithms and in all simulations by randomly sampling once from each set of feasible portfolios.

\subsubsection{Mutation effects on LLSO variants}

Our first task is to study the impact of the developed mutation operator on the LLSO-type algorithms, all equipped with our hybrid constraint-handling technique. For this purpose, we compare the following variants of the LLSO: the dynamic LLSO (DLLSO) \cite{Y-C-D-L-G-Z}, the adaptive LLSO (ALLSO) \cite{S-Y-G-M-L-Z} and the reinforcement learning level-based particle swarm optimization (RLLPSO) algorithm \cite{W-W-S}. The parameters setting of each algorithm is chosen following the literature, as reported in Table \ref{Table:2}.
\begin{table}[H]
\begin{center}
\caption{Parameter settings of the algorithms used in the comparisons.} \label{Table:2}
{\def\arraystretch{1.1}
{\small
\begin{tabular}{lccccccc}
\hline
Algorithm & Parameter settings & Reference \\
\hline
DLLSO & $S = \{4, 6, 8,10, 20, 50\}$, $NP = 500$, $\phi = 0.4$ & \cite{Y-C-D-L-G-Z}\\
ALLSO & $NL_{min} = 2$, $NL_{max} = 50$, $NP = 500$, $\bar{\delta} = 0.01$, & \\
& $px = 0.01$, $\xi = 10^{-6}$ & \cite{S-Y-G-M-L-Z} \\
RLLPSO & $S = \{4, 6, 8,10, 20, 50\}$, $NP = 500$, $\phi = 0.4$, & \\
& $\alpha = 0.4$, $\gamma = 0.8$, $\varepsilon = 0.9$, $\xi = 10^{-6}$ & \cite{W-W-S} \\
PSO & $\omega_{min} = 0.4$, $\omega_{max} = 0.9$, $c_{1, min} = c_{2, min} = 0.5$, & \\
& $c_{1, max} = c_{2, max} = 2.5$, & \cite{R-H-W},  \cite{C-T-F-P} \\
FA & $\alpha = 0.5$, $\beta_{min} = 0.2$, $\gamma = 1$ & \cite{Y} \\
\hline
\end{tabular}}
}
\end{center}
\end{table}

\begin{figure}[H]
\begin{center}
\includegraphics[width=.32\textwidth]{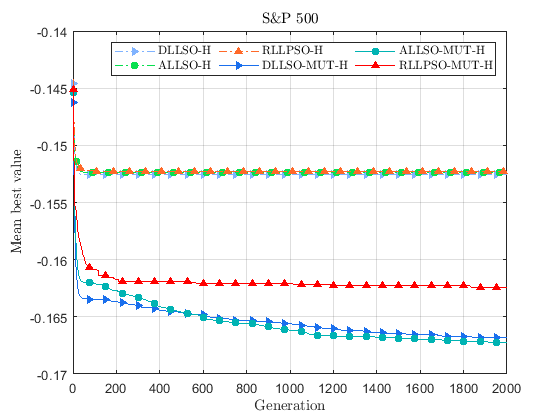}\includegraphics[width=.32\textwidth]{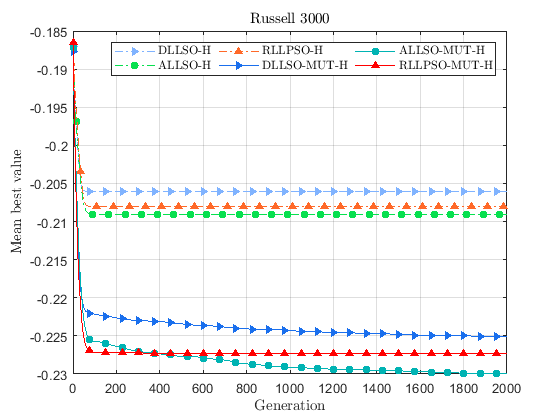}
\includegraphics[width=.32\textwidth]{Fig_comparison1_a3} \\
\includegraphics[width=.32\textwidth]{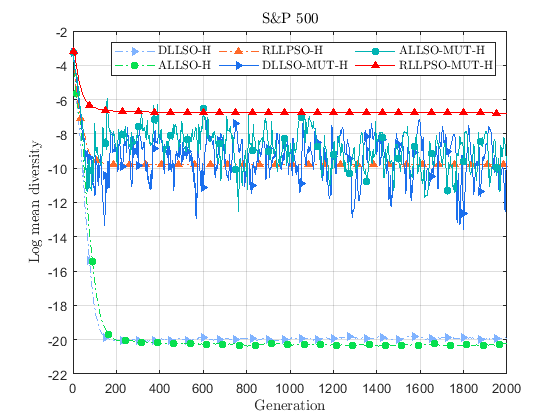}\includegraphics[width=.32\textwidth]{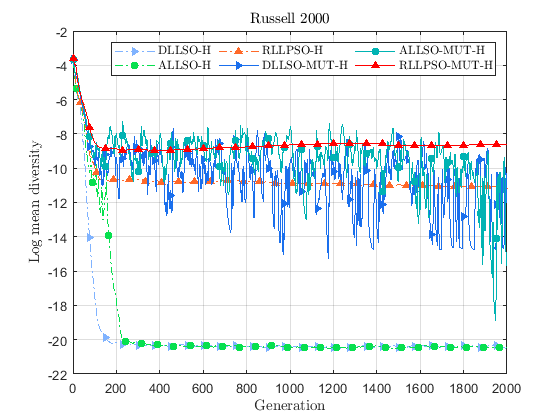}
\includegraphics[width=.32\textwidth]{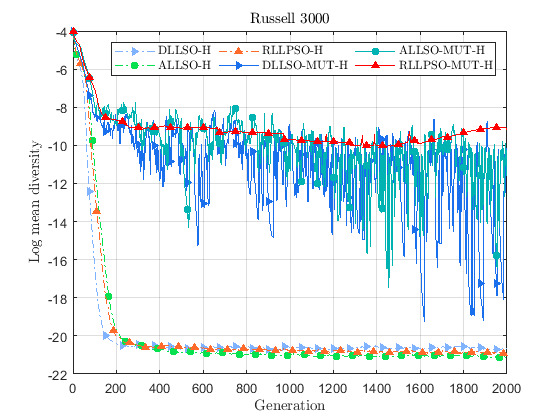} \\
\captionsetup{width=.9\linewidth}
\caption{Convergence and diversity analyses on the three data sets. Graphs in the first row show the behaviour of  algorithms in terms of mean best value of the objective function, while in the row below are displayed the logarithmic scale plots of the diversity scores.}
\label{FIG:OUT_1}
\end{center}
\end{figure}

Table \ref{tab:LLSO_comp} shows four performance metrics linked to the best objective function value over the 25 runs, on the three public data sets. The best results are highlighted in bold font.
We note that the ALLSO-MUT-H outperforms the competitors in all the case studies, presenting the lowest mean objective function value. Further, we remark that all the inspected LLSO variants are able to find feasible solutions.
Focusing on the mutation benefits, we observe from Table \ref{tab:improve} that the mutation has a significant impact on the performance of solvers. Specifically, although mutation-based optimizers present higher volatility than their counterparts, they always show lower results in terms of minimum-maximum range of the best solutions.
\begin{table}[H]
\begin{center}
\caption{Statistics regarding the best values of the objective function over the 25 runs.}\label{tab:LLSO_comp}
\resizebox{\textwidth}{!}{
{\def\arraystretch{1.3}
\resizebox{\textwidth}{!}{
\begin{tabular}{@{}clcccccc@{}}
\hline
\multicolumn{1}{c}{Data set} & Statistics & \multicolumn{1}{c}{DLLSO-H} & \multicolumn{1}{c}{DLLSO-MUT-H} & \multicolumn{1}{c}{ALLSO-H} & \multicolumn{1}{c}{ALLSO-MUT-H} & \multicolumn{1}{c}{RLLPSO-H} & \multicolumn{1}{c}{RLLPSO-MUT-H} \\
\hline
\multirow{ 4}{*}{\rotatebox{90}{S\&P 500}} & mean & -0.1525 & -0.1668 & -0.1524 & {\bf -0.1673} & -0.1523 & -0.1624 \\
& std & 0.0004  & 0.0015 & {\bf 0.0003} & 0.0013 & 0.0004 & 0.0024 \\
& min & -0.1532 & {\bf -0.1718} & -0.1529 & -0.1698 & -0.1532 & -0.1662 \\
& max & -0.1517 & {\bf -0.1643} & -0.1517 & -0.1638 & -0.1515 & -0.1546 \\
\hline
\multirow{ 4}{*}{\rotatebox{90}{Russell 2000}} & mean & -0.1921  & -0.2101 & -0.1925 & {\bf -0.2122} & -0.1929 & -0.2046 \\
& std & {\bf 0.0010} & 0.0043 & 0.0012 & 0.0037 & 0.0013 & 0.0040 \\
& min & -0.1932 & -0.2193  & -0.1946 & {\bf -0.2196} & -0.1963 & -0.2130 \\
& max & -0.1895 & -0.2003 & -0.1895 & {\bf -0.2039} & -0.1907 & -0.1959 \\
\hline
\multirow{ 4}{*}{\rotatebox{90}{Russell 3000}} & mean & -0.2060 & -0.2251 & -0.2090 & {\bf -0.2300} & -0.2080 & -0.2274 \\
& std & 0.0011 & 0.0032 & 0.0015 & 0.0025 & {\bf 0.0008} & 0.0043 \\
& min & -0.2085 & -0.2328 & -0.2120 & -0.2352 & -0.2094 & {\bf -0.2358} \\
& max & -0.2042 & -0.2192 & -0.2067 & {\bf -0.2260} & -0.2065 & -0.2196 \\
\hline
\end{tabular}
}}
}
\end{center}
\end{table}
\begin{table}[H]
\centering
\captionsetup{width=.9\linewidth}
\caption{Relative change of the mutated algorithms versus non-mutated counterparts. The $p$-values for the paired $t$-tests are displayed in brackets. Note that in all cases the p-values are under the significance level $\alpha=0.05$, indicating the rejection of the null hypothesis of equality of the means, against the alternative left-sided hypothesis.}
{\def\arraystretch{1.3}
\begin{tabular}{lccc}
\hline
& DLLSO-MUT-H & ALLSO-MUT-H & RLLPSO-MUT-H \\
Data set & vs. & vs. & vs. \\
& DLLSO-H (\%) & ALLSO-H (\%) & RLLPSO-H (\%) \\
\hline
S\&P 500 & 9.3705 & 9.7948 & 6.6287 \\
& ($6.8333 \cdot 10^{-25}$) & ($4.8584 \cdot 10^{-26}$) & ($7.5173 \cdot 10^{-17}$) \\
Russell 2000 & 9.3934 & 10.2666 & 6.0365 \\
& ($2.3999\cdot 10^{-17}$) & ($2.3049\cdot 10^{-19}$) & ($7.7679 \cdot 10^{-13}$) \\
Russell 3000 & 9.2573 & 10.0354 & 9.3067 \\
& ($2.0685\cdot 10^{-19}$) & ($1.913\cdot 10^{-23}$) & ($1.0664\cdot 10^{-17}$) \\
\hline
\end{tabular}}
\label{tab:improve}
\end{table}
Figure \ref{FIG:OUT_1} shows the convergence and the diversity analyses of the compared solvers on the three data sets.
From the first set of graphs, we note that the three mutated algorithms are able to reach significantly lower objective function values, and the ALLSO-MUT-H performs better than the others. Moreover, the algorithms without mutation show population stagnation around $100$ generations, meaning that they converge to a local minimum and are not able to further explore the search space.
This is confirmed by the results showed in the logarithmic scale plots of the diversity measures. We can observe that the ALLSO-MUT-H and the DLLSO-MUT-H are able to escape from the local minima, due to the oscillatory behaviour of the swarm diversity.

\subsubsection{Comparison with state-of-the-art swarm optimization algorithms}
In the previous subsection, we have analysed the impact of the mutation on the capabilities of LLSO-based algorithms, finding that the ALLSO-MUT-H is the more efficient choice in terms of convergence and quality of solutions. Now, the aim is to compare the above quoted optimizer with other state-of-the-art swarm optimization algorithms, namely the PSO \cite{C-T-F-P} and the Firefly algorithm (FA) \cite{Y}, both endowed with an exact $\ell_1$-penalty function. In this regard, the literature presents a wide range of penalty methods to tackle constraint-handling problems \cite{COE}. However, in the context of portfolio optimization allocation, the technique proposed in \cite{C-T-F-P} is the only one for which the convergence is guaranteed. For this reason, we have decided to adopt the exact $\ell_1$-penalty technique.

\begin{figure}[H]
\begin{center}
\includegraphics[width=.32\textwidth]{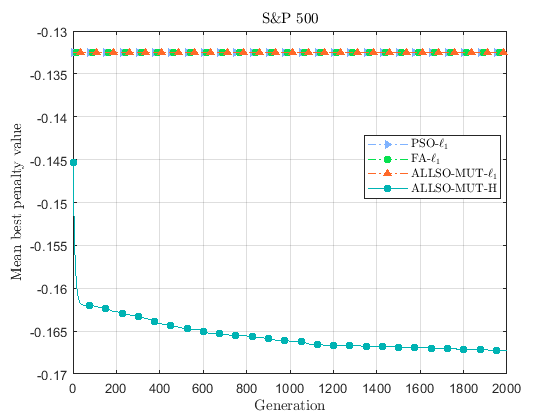}\includegraphics[width=.32\textwidth]{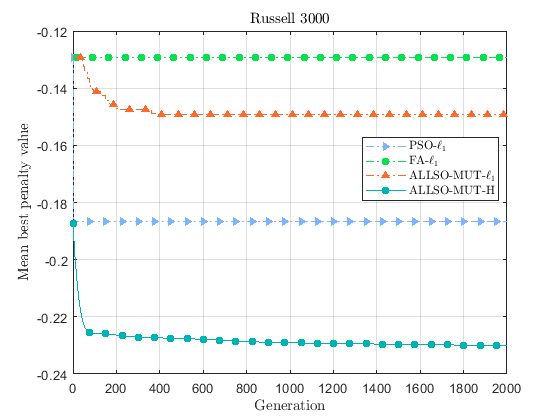}
\includegraphics[width=.32\textwidth]{Fig_comparison2_a3} \\
\includegraphics[width=.32\textwidth]{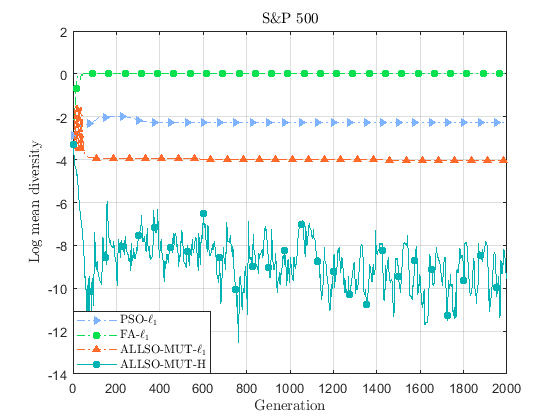}\includegraphics[width=.32\textwidth]{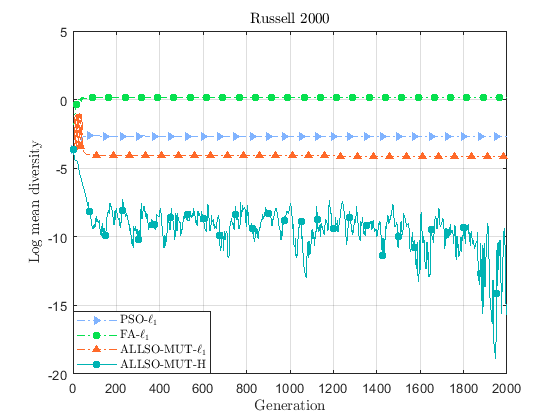}
\includegraphics[width=.32\textwidth]{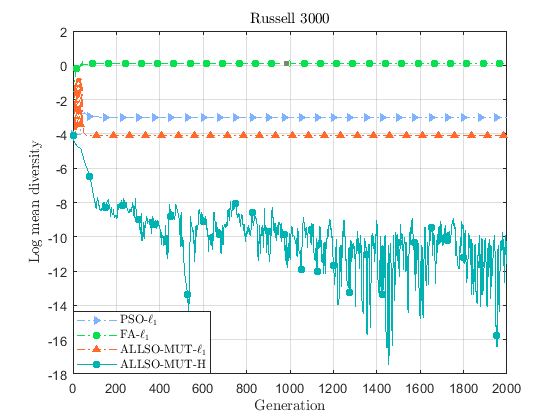} \\
\captionsetup{width=.9\linewidth}
\caption{Plots in the first row show the behaviour of the algorithms in terms of mean best value of the penalty function, while in the second row are presented the logarithmic scale graphs of the diversity.}
\label{FIG:OUT_2}
\end{center}
\end{figure}

We recall also that we adopt the same parameter setup presented above in Table \ref{Table:2}.
We exhibit the statistics of the comparison in Table \ref{tab:L1_comp}, where are displayed the percentage of feasible solutions provided by the different solvers over the 25 runs; the mean of the constraint violation function CV for the non-feasible solutions; the average value of the penalty function $F_{\ell_1}$ over the 25 runs. Notice that, the penalty function corresponds to the objective function when the solutions are feasible.
Looking at the results, we can argue that the ALLSO-MUT-H reaches the best mean value of the penalty function in all the data sets, and it always provides feasible solutions.
This insight is confirmed by the convergence analysis plots in Figure \ref{FIG:OUT_2}, which show the benefits of our constraint-handling technique in terms of accuracy of solutions. Moreover, the diversity graphs suggest that our solver is the sole algorithm able to exhibit exploration and exploitation phases alternatively.
\begin{table}[H]
\begin{center}
\captionsetup{width=.9\linewidth}
\caption{Comparison with state-of-the-art swarm optimization algorithms implementing the exact $\ell_1$-penalty framework proposed in \cite{C-T-F-P}.}\label{tab:L1_comp}
\resizebox{\textwidth}{!}{
{\def\arraystretch{1.6}
\begin{tabular}{@{}clcccc@{}}
\hline
\multicolumn{1}{c}{Data set} & Statistics & \multicolumn{1}{c}{PSO-$\ell_1$} & \multicolumn{1}{c}{FA-$\ell_1$} & \multicolumn{1}{c}{ALLSO-MUT-$\ell_1$} & \multicolumn{1}{c}{ALLSO-MUT-H} \\
\hline
\multirow{ 3}{*}{\rotatebox{90}{S\&P 500}} & feasible sol. (\%) & 0 & 0 & 100 & 100 \\
& mean CV & $3.1572\cdot 10^{-12}$ & $3.1572\cdot 10^{-12}$ & 0 & 0 \\
& mean $F_{\ell_1}$ & -0.1325 & -0.1325 & -0.1325 & -0.1673 \\
\hline
\multirow{ 3}{*}{\rotatebox{90}{Russell 2000}} & feasible sol. (\%) & 68 & 100 & 100 & 100 \\
& mean CV & $1.6035\cdot 10^{-5}$  & 0 & 0 & 0 \\
& mean $F_{\ell_1}$ & -0.1718 & -0.1639 & -0.1639 & -0.2122 \\
\hline
\multirow{ 3}{*}{\rotatebox{90}{Russell 3000}} & feasible sol. (\%) & 80 & 100 & 0 & 100 \\
& mean CV & $7.4015 \cdot 10^{-12}$  & 0 & 0.0320 & 0 \\
& mean $F_{\ell_1}$ & -0.1866 & -0.1291 & -0.1493 & -0.2300 \\
\hline
\end{tabular}
}
}
\end{center}
\end{table}

\subsection{Real-world application}

\subsubsection{Data description and investment setting}

The constituents of the MSCI World index, on 31st January 2022, form our investible pool. The data set has been downloaded from DataStream and consists of monthly prices covering the period from January 2012 to January 2022 for a total of 121 months. Stocks with missing observations were disregarded, and thus the final data set includes 1119 stocks. For the performance comparisons, we introduce a value-weighted benchmark index with the same constituents and, as in the previous analysis, we set $r_f = 0$.

The portfolio design employs the following parameter setting. For the cardinality constraint \eqref{card}, we consider $k \in \{335,167,55,22\}$, corresponding to $k_{\%} = 30\%$, $15\%$, $5\%$ and $2\%$ of the pool size, respectively. As stated in the introductory part of Section \ref{subsec:algo_comparisons}, we recall that the box thresholds in \eqref{box} are $l_i=0.001$ and $u_i=0.05$ for each asset $i$, with $i = 1,\ldots, 1119$, and the turnover rate in \eqref{turn} is set equal to $0.20$, as in \cite{K-B-C}.

We use a rolling time window procedure to rebalance optimal portfolios every month, from January 2017 to January 2022, to point out the effects of the market changes on the behaviour of the investments and, as a consequence, the total number of ex-post dates is 61. We solve the corresponding problem instances by employing overlapping 60-months windows, which are updated every month by removing the oldest data and including the latest information.

In each quoted window, as already pointed out above, we adopt a historical approach to calculate expected rates of return, and to reduce the bias in the estimation of the covariance matrix $C$, we take advantage of the shrinkage estimator proposed in \cite{L-W}. Let us denote by $\vec x_{t}$ the optimal portfolio at the ex-post month $t$, with $t = 1, \ldots, 61$. Due to the time dependence of the considered investment plan, we rewrite the turnover constraint \eqref{turn} as follows
\begin{equation}\label{turn_t}
\sum_{i=1}^{n} \lvert x_{t,i} - x_{t^-,i} \rvert \leq TR\, .
\end{equation}
In the previous equation, $\vec x_{t^-} = (x_{t^-,1}, \ldots, x_{t^-,n})$ represents the portfolio to be rebalanced \cite{S-W-M}, which is defined for $t = 2, \ldots, 61$ as
\begin{equation}\label{norm}
x_{t^- ,\, i} = \frac{x_{t-1 ,\, i} R^g_{t-1 ,\, i}}{\sum_{j = 1}^n x_{t-1 ,\, j} R^g_{t-1 ,\, j}}
\end{equation}
with the denominator being the gross portfolio return at month $t-1$.\footnote{The gross return of asset $i$ at month $t$ is defined as $R^g_{i,t} = \frac{S_{i,t}}{S_{i,t-1}}$, where $S_{i,t}$ is the price of the $i$-th asset at the end of month $t$.} At time $t = 1$, we set $\vec x_{t^-} = \vec 0$ and $TR = 1$.

Let us assume a self-financing strategy with an initial wealth $W_0 = 10,000,000$ \$. Then, we explicitly evaluate the magnitude of the trading through the cost function $\lambda(\vec x_t, \vec x_{t^-})$ introduced in \cite{B-V-F-P}. As reported in Table \ref{tab:cost}, we consider the transaction cost structure characterized by decreasing cost rates as the traded value increases.

\begin{table}[H]
\centering
\caption{Structure of transaction costs.}
\begin{tabular}{ccc}
\hline
Trading segment (\$) & Fixed fee (\$) & Proportional cost (\%)\\
\hline
0 -- 7,999 & 40 & 0\\
8,000 -- 49,999 & 0 & 0.5\\
50,000 -- 99,999 & 0 & 0.4\\
100,000 -- 199,999 & 0 & 0.25\\
$\geq 200,000$ & 400 & 0\\
\hline
\end{tabular}
\label{tab:cost}
\end{table}

\subsubsection{Ex-post performance measures}

The following measures are considered to evaluate the profitability of the investment strategies. Let $r_{p,t}^{out}$ be the ex-post portfolio rate of return realized at time $t$, with $t=1,\ldots, 61$. First, we consider the so-called ex-post Sharpe ratio \cite{Sh_1994}, defined as
\begin{equation}\label{sr}
SR^{out} = \frac{\mu^{out}}{\sigma^{out}}
\end{equation}
where $\mu^{out}$ and $\sigma^{out}$ are the mean and the standard deviation of the ex-post portfolio rates of return, respectively.\\
The second measure employed in the analysis is the so-called Omega ratio \cite{K-S}, defined as the ratio between the gains over a threshold level and the losses under a threshold level. In this study, we set both thresholds equal to zero, that is
\begin{equation}\label{omega}
Omega = \frac{\sum_{t=1}^{61}r_{p,t}^{out}\bbbone_{\{r_{p,t}^{out}>0\}}}{-\sum_{t=1}^{61}r_{p,t}^{out}\bbbone_{\{r_{p,t}^{out}<0\}}}
\end{equation}
where $\bbbone_{A}$ is the indicator function on $A$.\\The information gathered from these performance measures draws a complete picture of the ex-post portfolio return distribution. In particular, the ex-post Sharpe ratio describes the central part of the portfolio return distribution, while the Omega ratio considers the behaviour of profits and losses.

Further, to measure the profitability of the investment at time $t$, we compute the net wealth as
\begin{equation}\label{wealth}
W_{t} = W_{t-1}\left (1+r_{p,t}^{out} \right ) - \lambda(\vec x_t, \vec x_{t^-}).
\end{equation}
Then, we compare the profitability of the investments using the so-called compound annual growth rate \cite{Sh_1994}, which in our case is calculated as
\begin{equation}\label{cagr}
CAGR = \left(\frac{W_{61}}{W_0}\right)^{\frac{12}{61}}-1
\end{equation}
where $W_0$ represents the initial wealth and $W_{61}$ is the final wealth.\\
To evaluate the capacity of a strategy to avoid high losses, we introduce the drawdown measure \cite{C-S-U-Z}, which can be written
\begin{equation}\label{dd}
DD_t = \min \left \{0,\frac{W_t - W_{peak}}{W_{peak}} \right \}
\end{equation}
where $W_{peak}$ is the maximum amount of wealth reached by the strategy until time $t$. In particular, we consider the mean and the standard deviation of the drawdown measure over time.

Finally, we propose to measure the effect of the costs on the available capital in the out-of-sample period by
\begin{equation}\label{EQ:costs1}
\Lambda_{\%} = \frac{1}{61}\sum_{t = 1}^{61} \frac{\lambda(\vec x_t, \vec x_{t^-})}{W_{t-1}} \cdot 100.
\end{equation}

\subsubsection{Ex-post performance analysis}

In the ex-post analysis, we investigate how the performance of the proposed asset allocation model changes by varying the $k$ parameter in \eqref{card}.

First, we remark that, for any ex-post dates, the proposed hybrid LLSO variant identifies feasible solutions for all the portfolio sizes.
The empirical results are summarised in Table \ref{Table:RealWorld}, where the number of assets of the considered strategies is also displayed. Note that for each value of $k$, the proposed investments provide better performances than the value-weighted benchmark. This result implies that introducing a cardinality constraint in the portfolio model allows to choose a subset of the most profitable assets in the investible pool.

In terms of the return-risk profile, strategies with $k_{\%} = 30\%$, $15\%$, and $5\%$ show comparable performances, while the strategy with $k_{\%} = 2\%$ has a lower Sharpe ratio, which is due to its large volatility. Similar conclusions can be made about the Omega ratio, which expresses the gain-loss profile of the strategies. Despite better performance with respect to Sharpe and Omega ratios, strategies involving portfolios with a larger number of assets generate less wealth. Moreover, we observe that reducing cardinality leads to more profitable portfolio strategies.

Concerning the drawdown measures, the $5\%$ asset allocation model is the most conservative, while the one with $k_{\%} = 2\%$ is the worst. Thus, the performance deteriorates by reducing portfolio size below a critical threshold.

\begin{figure}[H]
\begin{center}
\includegraphics[width=0.8\textwidth]{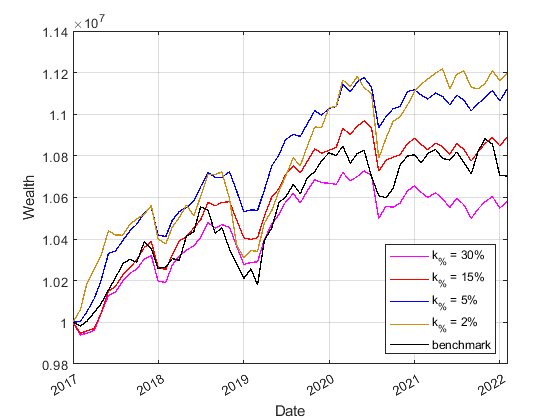}
\captionsetup{width=.9\linewidth}
\caption{Ex-post evolution of net wealth of the benchmark and of the proposed cardinality-constrained portfolio allocation model with different cardinalities.}
\label{FIG:OUT_wealth}
\end{center}
\end{figure}

As highlighted in the last two rows of Table \ref{Table:RealWorld} and in Figure \ref{FIG:OUT_wealth}, the impact of transaction costs for strategies with small $k$ is negligible. On the contrary, portfolios with many assets have more fluctuations in the rebalancing phases, leading to higher trading commissions with a significant impact on the wealth generated.

Summing up, we can infer that the strategy with $k_{\%} = 5\%$ shows the best balance between risk-adjusted performance measures and capability to generate net profits.

\begin{table}[H]
\begin{center}
\captionsetup{width=.9\linewidth}
\caption{Performance of the proposed cardinality-constrained portfolio allocation model for different cardinalities in comparison to the benchmark.}
\label{table:Performance:out-of-sample1} \label{Table:RealWorld}
{\small
\begin{tabular}{lccccccc}
\hline
$k_{\%}$ & 30\% & 15\% & 5\% & 2\% & Benchmark \\
\hline
num. assets & 335 & 167 & 55 & 22 & 1119 \\
Sharpe ratio & 0.4917 & 0.4794 & 0.4467 & 0.2878 & 0.1731 \\
Omega ratio & 3.3269 & 3.2428 & 3.0360 & 2.0980 & 1.5388 \\
CAGR & 1.1148 & 1.6899 & 2.1118 & 2.2502 & 1.3460 \\
std & 0.0052 & 0.0052 & 0.0054 & 0.0077 & 0.0066 \\
mean DD & -0.0066 & -0.0057 & -0.0050 & -0.0072 & -0.0070 \\
std DD & 0.0070 & 0.0064 & 0.0060 & 0.0106 & 0.0089 \\
mean $\lambda$ (\$) & 16,046 & 10,713 & 6,571.7 & 3,121.9 & -- \\
$\Lambda_{\%}$ & 0.1536 & 0.1011 & 0.0610 & 0.0290 & -- \\
\hline
\end{tabular}
}
\end{center}
\end{table}

\section{Conclusions and future works}\label{sec:conc}
In this paper, we have developed a swarm optimization algorithm for solving a large-scale cardinality-constrained portfolio optimization problem, where the modified Sharpe ratio performance measure represents the objective function. We have considered four real-world constraints: cardinality, box, budget and turnover constraints. Due to the properties of the model inspected,
we have proposed a variant of the LLSO equipped with a hybrid procedure to manage the constraints efficiently. Moreover, a novel mutation operator has been introduced to improve the accuracy of solutions. Our solver capabilities have been compared with those of two variants of the LLSO as well as other state-of-the-art swarm optimization algorithms endowed with an $\ell_1$-penalty function. Numerical experiments on three publicly available large-scale data sets showed the outperformance of our hybrid procedure. From the financial point of view, we have analysed the sensitivity of the portfolio model to the cardinality constraint with data from the last five years of the MSCI World index. We have found that portfolios of small size are more competitive with respect to the value-weighted benchmark index, also in periods of market downturns. Specifically, the losses are reduced, and the cost impact on the available capital is marginal compared to the profits.
In the future works, we plan to further assess the capabilities of the proposed hybrid constraint-handling technique by including other constraints in the asset allocation model. On the one hand, we will consider the so-called risk-budgeting constraints to explicitly control the portfolio risk exposition. On the other hand, based on the European Green Deal and the ESG Disclosure requirements for funds and investments, we will add sustainable-policy constraints to guarantee a minimum level of ESG rating to the investment. In addition, we will consider alternative performance measures in the optimization framework to properly handle tail-risk, and we will extend the experimental part by using different data sets.


\newpage

\appendix
\section{Appendix -- Proofs of the main results} \label{Appendix_A}

\begin{proof}{\ref{prop_proj}}
We rewrite problem \eqref{probl_components} in the following way
\begin{equation*}
\min_{\vec x_K: \, x_i\in C_i}\left\{\frac{1}{2}\sum_{j\in K}(x_j-y_j)^2\right\}+\frac{1}{2}\sum_{j\in I\setminus K}y_j^2\, ,
\end{equation*}
that is equivalent to
\begin{equation}\label{prob_min_sum}
\min_{\vec x_K: \, x_i\in C_i}\left\{\frac{1}{2}\sum_{j\in  K}(x_j-y_j)^2\right\}-\frac{1}{2}\sum_{j\in  K}y_j^2+\frac{1}{2}\sum_{j\in I}y_j^2\, .
\end{equation}
We note that the last term in \eqref{prob_min_sum} does not depend on $x_j, \, j\in K$, so we can focus our attention on the first two terms, i.e.
$$ \min_{\pi_K(\vec x): \, x_i\in C_i}\frac{1}{2}\|\pi_K(\vec x - \vec y)\|^2-\frac{1}{2}\|\pi_K(\vec y)\|^2\, . $$
By contradiction, we suppose that there is a $K^\prime$ different from $ K^\ast$, where we recall that $K^\ast$ is the set of indices corresponding to the $k$ largest components of $\vec y$.
At this point, we define
\begin{equation*}
f(\pi_K(\vec y))= -\frac{1}{2}\|\pi_K(\vec y)\|^2+\min_{\pi_K(\vec x):\, x_i\in C_i}\frac{1}{2}\|\pi_K(\vec x - \vec y)\|^2
\end{equation*}
and
\begin{equation*}
g(t)= f((1-t)\pi_{K^\ast}(\vec y)+t\pi_{K^\prime}(\vec y))\quad \text{ with }\quad t\in [0,1]\, .
\end{equation*}
Then we have,
\begin{equation*}
f(\pi_{K^\prime}(\vec y))-f(\pi_{K^\ast}(\vec y))=g(1)-g(0)=\int_0^1 g^\prime(t) \dt
\end{equation*}
and
\begin{equation*}
g^\prime(t)=\nabla f((1-t)\pi_{K^\ast}(\vec y)+t\pi_{K^\prime}(\vec y))\cdot (-\pi_{K^\ast}(\vec y)+\pi_{K^\prime}(\vec y))\, ,
\end{equation*}
where $\nabla f(\pi_{K}(\vec y))=-\pi_{K}(\vec y)+\pi_{K}(\vec y)-\pi_{K}(\vec x^\ast) =  -\pi_{K}(\vec x^\ast) $ and $\pi_{K}(\vec x^\ast)=\displaystyle \argmin _{\pi_{K}(\vec x):\, x_i\in C_i}\frac{1}{2}\|\pi_{K}(\vec x - \vec y)\|^2$.
\medskip
Now, since $\mc C_i\subset \R_+$, we have that $\nabla f(\pi_{K}(\vec y))$ is non-positive in all components. Moreover, $-\pi_{K^\ast}(\vec y)+\pi_{K^\prime}(\vec y)\leq 0$ due to the fact that $\pi_{K^\ast}(\vec y)$ is the projection of $\vec y$ onto the set of its largest components. As a result, we obtain $g^\prime (t)\geq 0 $, which
implies $f(\pi_{K^\prime}(\vec y))\geq f(\pi_{K^\ast}(\vec y))$. This means that $ K^\ast$ must be the optimal choice.
\qed
\end{proof}

\bigskip

\begin{proof}{\ref{prop_equiv}}
The proof of the first part of the proposition follows by defining $\vec x^{\ast}_{K^\ast}$ such that $\vec x^{\ast}_{K^\ast} = \vec \delta^\ast \otimes \vec x^\ast$, where $\otimes$ stands for the Hadamard product.\\On the contrary, if $\vec x^{\ast}$ solves \eqref{opt_bis}, then $\vec x^{\ast}_{K^\ast} = \vec x^{\ast}$. By taking
\begin{equation*}
\delta_i^\ast =\begin{cases}
1\quad \text{if} \quad i\in K^\ast\\
0 \quad \text{otherwise}
\end{cases}
\end{equation*}
we deduce that $\left(\vec \delta^\ast ,\,\ \vec x^\ast\right)$ solves \eqref{opt1}.\qed
\end{proof}

\bigskip

\begin{proof}{\ref{prop:proj_beck}}
The result follows from Theorem 6.27 in \cite{B}, where the hyperplane is represented by the budget constrain \eqref{budg} and the box is $[ \pi_{K^{\ast\ast}}(\vec l) ,\,  \pi_{K^{\ast\ast}}(\vec u)]$. We recall also that $\pi_{K^{\ast\ast}}(\vec z)= (z_{i_1}, \ldots ,z_{i_k})$ with $z_{i_j}>0$ and $i_j\in K^{\ast \ast}$. \qed
\end{proof}

\section{Appendix -- An approach based on the exact $\ell_1$-penalty function} \label{Appendix_B}
In this appendix we introduce a procedure based on the exact $\ell_1$-penalty function for solving cardinality-constrained portfolio optimization problems. We adapt the approach discussed in \cite{C-T-F-P} to the algorithms used in the comparison analysis of Subsection \ref{subsec:algo_comparisons} and thus we define the constraint violations as follows
\begin{equation*}
\begin{split}
CV_1& = \left|\sum_{i=1}^n x_i-1\right|\\
CV_2&= \max \left\{ \sum_{i=1}^n \delta_i -k, 0 \right\}\\
CV_3&=\sum_{i=1}^n \max \left\{ \delta_i l_i -x_i, 0 \right\}\\
CV_4&= \sum_{i=1}^n \max \left\{ x_i - \delta_i u_i , 0 \right\}\\
CV_5&=\sum_{i=1}^n \left|\delta_i (1-\delta_i)\right|\\
CV_6&= \max\left\{\sum_{i=1}^n\left|x_i-x_{0,i}\right|-TR,\, 0\right\} \, .
\end{split}
\end{equation*}
In this manner, we introduce the exact $\ell_1$-penalty function
\begin{equation}\label{l1penalty}
F_{\ell_1}(\vec{x} , \vec{\delta}; \vec{\veps}) = f(\vec{x})+ \frac{1}{\veps_0}\left[\veps_1 CV_1+\veps_2 CV_2
+\veps_3 CV_3+\veps_4 CV_4
+\veps_5 CV_5+\veps_6 CV_6\right]
\end{equation}
where $\vec{\veps}=\left(\veps_0, \veps_1, \ldots ,\veps_6\right)$, with $\veps >0 $ for all $i$.

The initial parameters vector $\vec{\veps}^0$ is set to
$\vec{\veps}^0=\left(\veps^0_0, \veps^0_1, \ldots ,\veps^0_6\right)=(10^{-4},1,\ldots ,1)\in \R^7$,
where $\veps_0^0$ is chosen in order to privilege feasible solutions, and the other parameters are equally penalized for all constraint violations.

The vector $\vec \veps$ is updated by checking the decrease of the function $f(\vec{x})$ and the violation of the constraints. More precisely, on the one hand, every $5$ iterations the entry $\veps_0(g)$ is updated according to the rule
\begin{equation} \label{eqn:eps_0}
\veps_0(g+1) = \begin{cases}
\min \{3 \cdot \veps_0(g), \, 1 \} \quad \text{ if }\quad f(\vec{x}(g))\geq f(\vec{x}(g-1))\\
\max \{0.6 \cdot \veps_0(g), \, 10^{-15} \} \quad \text{ if }\quad f(\vec{x}(g))< 0.9 \cdot f(\vec{x}(g-1))\\
\veps_0(g) \quad \text{ otherwise}.
\end{cases}
\end{equation}
On the other hand, every $10$ iterations the entries $\veps_i(g)$, $i=1, \ldots 6$, are updated following the scheme
\begin{equation} \label{eqn:eps_i}
\veps_i(g+1)=\begin{cases}
\min \{2 \cdot \veps_i(g), \, 10^4 \} \quad \text{ if }\quad CV_i(g)> 0.95\cdot CV_i(g-1)\\
\max \{0.5 \cdot \veps_i^g, \, 10^{-4} \} \quad \text{ if }\quad CV_i(g) < 0.9 \cdot CV_i(g-1) \\
\veps_i(g) \quad \text{ otherwise},
\end{cases}
\end{equation}
with $CV_i$ the respective constraint violation linked to the $\veps_i$ parameter.

The above quoted strategy privileges
optimality of solutions possibly at the expenses of their feasibility, due to the fact that $\veps_0(g+1)$ in \eqref{eqn:eps_0} is increasing in $F_{\ell_1}(\vec{x},\vec{\delta}; \vec{\veps}(g+1))$ when the function
value $f(\vec{x}(g))$ increases. Moreover, to favour feasibility of solutions possibly at the expenses of their optimality, the penalty parameter $\veps_i(g+1)$ in \eqref{eqn:eps_i} is increased when the relative constraint violation in the $g$-th generation increases with respect to the previous one.

The procedure is also equipped by a splitting and refining technique for the positions of the particles. In particular, at each iteration, a particle $p$ is split in its components $\vec{x}^p(g)$ and $\vec{\delta}^p(g)$ that are updated separately. For the vector $\vec{\delta}^p(g)$  we employ the following updating rule
\begin{equation*}
\delta_i^p(g+1)=\begin{cases}
1 \quad \text{ if }\quad x_i^p(g) \in [l_i,\, u_i]\\
0 \quad \text{ otherwise}
\end{cases}
\end{equation*}
for $i=1, \ldots , n$. Then, $\vec \delta^p(g+1)$ is kept fixed and $F_{\ell_1}(\vec{x}, \vec{\delta}^p(g+1);\vec \veps(g+1))$ is minimized with respect to $\vec{x}$, obtaining $\widetilde{\vec{x}}^p(g+1)$. Finally, $\widetilde{\vec{x}}^p(g+1)$ is refined getting
\begin{equation}
x_i^p(g+1) = \frac{\widetilde{x}_i^p(g+1) \delta_i^p(g+1)}{\sum_{i=1}^n \widetilde{x}_i^p(g+1) \delta_i^p(g+1)}\, ,
\end{equation}
for $j=1, \ldots , NP$.

\subsection*{Acknowledgements}
{\small

The third author is member of the INdAM (Italian Institute for Advanced Mathematics) group.
}




\end{document}